\documentclass[11pt]{article}
\usepackage[left=1in,top=1in,right=1in,bottom=1in,head=.1in,nofoot]{geometry}

\usepackage{setspace,url,bm,amsmath} 

\usepackage{titlesec} 
\titlelabel{\thetitle.\quad} 

\usepackage{graphicx} 
\usepackage{bbm}
\usepackage{latexsym}
\usepackage{caption}
\usepackage[margin=20pt]{subcaption}
\usepackage{hyperref}
\usepackage{enumerate}

\usepackage[table]{xcolor}

\newcommand{\GG}[1]{}

\usepackage{amsthm}
\usepackage{comment}
\theoremstyle{definition}

\newtheorem*{theorem*}{Theorem}

\newtheorem*{corollary*}{Corollary}

\usepackage{natbib} 
\bibpunct{(}{)}{;}{a}{}{,} 

\usepackage{etoolbox} 
\apptocmd{\sloppy}{\hbadness 10000\relax}{}{} 

\usepackage{color}
\usepackage{listings}

\def\ind{\begin{picture}(9,8)
         \put(0,0){\line(1,0){9}}
         \put(3,0){\line(0,1){8}}
         \put(6,0){\line(0,1){8}}
         \end{picture}
        }

\def\Expo{\text{Expo}}
\def\Laplace{\text{Laplace}}
\def\G{\text{Gamma}}
\def\d{\textnormal{d}}

\begin{document}
\doublespacing
\title{\bf Representation for the Gauss--Laplace Transmutation}
\author{Peng Ding and Joseph K. Blitzstein 
\footnote{
Peng Ding (Email: \texttt{pengdingpku@gmail.com}) is Assistant Professor at University of California at Berkeley, and Joseph K. Blitzstein (Email: \texttt{blitz@fas.harvard.edu}) is Professor of Practice of Statistics at Harvard University. 
The authors thank Dr. Avi Feller at Goldman School of Public Policy at Berkeley for bringing Professor Christian Robert's blog article on Gauss--Laplace transmutation (\url{https://xianblog.wordpress.com/2015/10/14/gauss-to-laplace-transmutation/}) into their attention. 
}}

\date{}
\maketitle

\begin{abstract}
Under certain conditions, a symmetric unimodal continuous random variable $\xi$ can be represented as a scale mixture of the standard Normal distribution $Z$, i.e., $\xi = \sqrt{W} Z$, where the mixing distribution $W$ is independent of $Z.$ It is well known that if the mixing distribution is inverse Gamma, then $\xi$ is student's $t$ distribution. However, it is less well known that if the mixing distribution is Gamma, then $\xi$ is a Laplace distribution. Several existing proofs of the latter result rely on complex calculus and change of variables in integrals. We offer two simple and intuitive proofs based on representation and moment generating functions.

\end{abstract}

\section{Main Result to Prove}

Let $\Expo(1)$ denote the standard Exponential distribution with mean one, $\Laplace$ the standard Laplace distribution with mean zero and variance two, $N(\mu,\sigma^2)$ the Normal distribution with mean $\mu$ and variance $\sigma^2.$ If two random variables $A$ and $B$ have the same distribution, then we write $A\sim B.$
The Gauss--Laplace transmutation states that
$$
V\sim 2\Expo(1),\quad L\mid V\sim N(0, V) \Longrightarrow L\sim \Laplace,
$$ 
or equivalently, if $\Expo\sim \Expo(1)$ is independent of $Z\sim N(0,1)$, then
$$
L = \sqrt{2 \Expo} Z \sim \Laplace.
$$
Interestingly, the Laplace and Normal distributions are Pierre-Simon Laplace's first and second law of errors \citep{wilson1923first}, which are tied together by a scaling factor $\sqrt{2 \Expo}.$ In high dimensional statistics, this result is crucial to efficiently simulate posterior distribution of the Bayesian Lasso \citep{park2008bayesian}, which imposes Laplace priors on the regression coefficients. Some proofs of the Gauss--Laplace transmutation exist in the literature \citep{andrews1974scale, west1987scale}, which rely on advanced theory of complex calculus. In this note, we offer two simple proofs based on representation and moment generating functions (MGFs).

\section{Review of Some Basic Representations}

Let $Z\sim N(0,1), \Expo\sim \Expo(1), X^2_k\sim \chi^2_k$, and $L\sim \Laplace$. We use these symbols for Normal, Exponential, Chi-Squared, and Laplace random variables from now on. We also use $S$ for a random sign taking values $\pm 1$ with probabilities $1/2.$
Let IID denote ``independently and identically distributed.'' The following representations are useful.

\begin{enumerate}
[(a)]
\item\label{rep:a}
If $Z_1, Z_2$ are IID $N(0,1)$, then $ 2\Expo\sim X^2_{2}\sim Z_1^2+Z_2^2$.

\begin{proof}
Let $\G(a)$ be a Gamma random variable with shape parameter $a$ and rate parameter one. 
We have $2\Expo\sim 2 \G(1) \sim X^2_2$ by the relationship between Gamma and Chi-Squared random variables, and $2\Expo\sim Z_1^2 + Z_2^2$ by the definition the Chi-Squared random variable. 
\end{proof}

\item\label{rep:b}
If $Z_1, Z_2$ are IID $N(0,1)$, then $Z_1Z_2\sim (Z_1^2-Z_2^2)/2.$

\begin{proof}
Because the bivariate Normal distributions satisfy $(Z_1,Z_2)\sim (Z_1+Z_2, Z_1-Z_2)/\sqrt{2}$, we have $Z_1Z_2\sim (Z_1+Z_2)( Z_1-Z_2)/2 = (Z_1^2-Z_2^2)/2.$
\end{proof}

\item\label{rep:c}
If $\Expo_1, \Expo_2$ are IID $\Expo(1)$ independent of $S$, then $L\sim S \Expo_1 \sim  \Expo_1 - \Expo_2.$

\begin{proof}
Obviously, $L\sim S \Expo_1$, which justifies the other name, double exponential, of the Laplace distribution. By symmetry and the memoryless property, we have 
$$
\Expo_1 - \Expo_2\sim S \{  \max(\Expo_1, \Expo_2)  - \min(\Expo_1, \Expo_2)   \} 
\sim  S  \Expo_1.
$$
\end{proof}

\end{enumerate}

\section{A Proof Based on Representation}

\begin{proof}
We let $Z, Z_1,Z_2,Z_3,Z_4$ be IID $N(0,1)$, and $\Expo_1,\Expo_2$ be IID $\Expo(1).$
Using representation (\ref{rep:a}), we have 
$$
L = \sqrt{2 \Expo} Z\sim \sqrt{ Z_1^2 + Z_2^2 } Z.
$$ 
By conditioning on $(Z_1,Z_2)$, we have
$$
Z_1Z_3 + Z_2Z_4  \sim N(0,    Z_1^2 + Z_2^2) \sim \sqrt{ Z_1^2 + Z_2^2 } Z,
$$
which implies 
$
L\sim Z_1Z_3 + Z_2Z_4.
$ 
Using representation (\ref{rep:b}), we have 
$$
L\sim (Z_1^2-Z_3^2)/2 + (Z_2^2-Z_4^2)/2 = (Z_1^2+Z_2^2)/2  - (Z_3^2+Z_4^2)/2 .
$$
Using representation (\ref{rep:a}) again, we have
$
L\sim \Expo_1 - \Expo_2.
$
Therefore, the final result follows from representation (\ref{rep:c}).
\end{proof}

\section{A Proof Based on Moment Generating Functions}

The MGF of $Z$ is $M_N(t) = E(e^{tZ}) = e^{t^2/2}$, and the MGF of $\Expo$ is $M_E(t) = E(e^{t\Expo}) = 1/(1-t)$ for $ t <1.$ We can prove the result by calculating the MGFs directly.

\begin{proof}
Using the known MGFs above and two conditional arguments, we can verify that
the MGF of $\sqrt{2 \Expo} Z $ is
$$
E( e^{t\sqrt{2 \Expo} Z } ) 
=E\{   E( e^{t\sqrt{2 \Expo} Z } \mid \Expo)  \}
=       E\{   M_N(t\sqrt{2 \Expo}) \}
=       E\{      e^{t^2 \Expo} \}
= M_E(t^2)
=     {1\over 1-t^2} ,
$$
and the MGF of $L\sim \Laplace\sim S \Expo$ is
$$
E(e^{tS\Expo}) 
= E\{  E(e^{tS\Expo} \mid S)  \}
=  E\{  M_E(tS) \}
=  E\left  ( 1 \over 1-tS \right) 
= \frac{1}{2}\left(  \frac{1}{1-t} + \frac{1}{1+t}  \right)
= {1\over 1-t^2} .
$$
\end{proof}

\section{Discussion}


Because of representation (\ref{rep:a}), we have $\sqrt{2 \Expo} Z \sim S \sqrt{X^2_1 X^2_2}$ and $L\sim S\Expo\sim S X_2^2/2$. Therefore, the Gauss--Laplace transmutation reduces to $4 X^2_1 X^2_2 \sim (X_2^2)^2$. In fact, a general result holds: $4 X^2_{k}  X^2_{k+1} \sim (X^2_{2k})^2.$ This is the stochastic analogue of Legendre's duplication formula for Gamma functions as discussed in \citet{gordon1989bounds}, \citet{gordon1994stochastic} and the Supplementary Materials.

The proof by representation also establishes the result $    \sqrt{X_2^2} Z\sim  Z_1Z_3 + Z_2Z_4$. In fact, a general result holds: $\sqrt{X_n^2} Z \sim \sum_{i=1}^n Z_{1i} Z_{2i}$, where $Z_{1i}, Z_{2i}\ (i=1,\ldots,n)$ are IID $N(0,1).$ 
This result is related to the Bartlett decomposition of the Wishart matrix \citep{anderson2003introduction} as commented in the Supplementary Materials.

As a byproduct of our proof by representation, we have shown that $ Z_1Z_3 + Z_2Z_4$ follows a Laplace distribution. Because of symmetry $Z_2\sim -Z_2$, the determinant of a two by two matrix with IID $N(0,1)$ entries, $Z_1Z_3 - Z_2Z_4$, also follows a Laplace distribution. This result has been established early in the literature \citep{nyquist1954distribution, nicholson1958distribution}, with a proof similar to our representation by \citet{mantel1966light}, and a proof based on characteristic functions by \citet{mantel1973characteristic}. In a series of exchanges in {\it The American Statistician}, these proofs were revived by \citet{missiakoulis1985distribution}, \citet{farebrother1986pitman} and \citet{mantel1987laplace}.

Professor Christian P. Robert posted online a proof based on directly calculating the density function. He generously allowed us to include his proof, which is in the Supplementary Materials. In a blog article about Gauss--Laplace transmutation, he was enthusiastic about proofs without complicated integrals or advanced theory, which motivated our note here.



%
%
%
%
%
%
%
%
%
%
%
%
%

\bibliographystyle{apalike}
\bibliography{GaussLaplace}

\newpage
\setcounter{page}{1}
\begin{center}
\bf \huge 
Supplementary Materials
\end{center}


\setcounter{equation}{0}
\setcounter{section}{0}
\setcounter{figure}{0}
\setcounter{example}{0}
\setcounter{proposition}{0}
\setcounter{theorem}{0}
\setcounter{table}{0}

\renewcommand {\theproposition} {A.\arabic{proposition}}
\renewcommand {\theexample} {A.\arabic{example}}
\renewcommand {\thefigure} {A.\arabic{figure}}
\renewcommand {\thetable} {A.\arabic{table}}
\renewcommand {\theequation} {A.\arabic{equation}}
\renewcommand {\thelemma} {A.\arabic{lemma}}
\renewcommand {\thesection} {A.\arabic{section}}
\renewcommand {\thetheorem} {A.\arabic{theorem}}

\section{Stochastic Analogue of Legendre's Duplication Formula}
If $X^2_k \ind X^2_{k+1}$, then $4 X^2_{k}  X^2_{k+1} \sim (X^2_{2k})^2.$

\begin{proof}
First, the Chi-Squared random variable has moment
$$
E(X^2_k)^t = 2^t \frac{\Gamma(k/2+t)}{\Gamma(k/2)} .
$$ 
Second, Legendre's Duplication Formula gives
$$
\Gamma(z)\Gamma(z+1/2) = 2^{1-2z} \sqrt{\pi} \Gamma(2z).
$$
We need only to show that $\log 4+\log X^2_{k}  + \log X^2_{k+1} $ and $2 \log X^2_{k+1}$ have the same MGFs. 

By definition, the MGF of the former is
$$
E\{  e^{ t (\log 4+\log X^2_{k}  + \log X^2_{k+1}  ) } \} = 4^t E (X^2_{k})^t E(X^2_{k+1})^t = 
4^t 2^t  \frac{\Gamma(k/2+t)}{\Gamma(k/2)}      2^t  \frac{\Gamma(k/2+t + 1/2)}{\Gamma(k/2 + 1/2)}
$$
and the MGF of the latter is
$$
E\{  e^{t(2 \log X^2_{2k})} \} = 2^{2t} E(  X_{2k}^2 )^{2t} = 2^{2t} 2^t \frac{\Gamma(k+t)}{\Gamma(k)}.
$$
In fact, the above two MGFs are the same according to Legendre's Duplication Formula.
\end{proof}

\section{Bartlett Decomposition of the Wishart Matrix}
 
If $ X_n^2\ind Z$, then $\sqrt{X_n^2} Z \sim \sum_{i=1}^n Z_{1i} Z_{2i}$, where $Z_{1i}, Z_{2i}\ (i=1,\ldots,n)$ are IID $N(0,1).$

This result is closely related to the Bartlett decomposition of the Wishart matrix. 
The sample covariance matrix of the $Z$'s, Wishart$_2(I, n)$, can be decomposed as
$$
\begin{pmatrix}
\sum_{i=1}^n Z_{1i}^2& \sum_{i=1}^n Z_{1i} Z_{2i} \\
\sum_{i=1}^n Z_{2i} Z_{1i} & \sum_{i=1}^n Z_{2i}^2
\end{pmatrix}
\sim 
 \begin{pmatrix}
 \sqrt{X_n^2} & 0\\
 Z & \sqrt{X_{n-1}^2}
 \end{pmatrix}
 \begin{pmatrix}
 \sqrt{X_n^2} & Z\\
 0 & \sqrt{X_{n-1}^2}
 \end{pmatrix} 
 =
 \begin{pmatrix}
 X_n^2 & \sqrt{X_n^2} Z\\
 \sqrt{X_n^2} Z & Z^2 +  X_{n-1}^2
 \end{pmatrix} . 
$$
Comparing the off diagonal term, we obtain the result.

\section{A Proof Based on Directly Calculating the Density Function}

The authors thank Professor Christian Robert for allowing us to include his proof posted on the website:
\url{http://stats.stackexchange.com/questions/175458}

With scale parameter $\lambda$, the Gauss--Laplace transmutation states that
$$
V\sim 2/\lambda^2\times \Expo(1),\quad L\mid V\sim N(0, V) \Longrightarrow L\sim  \lambda\times  \Laplace.
$$
Therefore, we need to verify the following integral:
$$
f(l) = \int_0^\infty \frac{1}{\sqrt{2\pi v} } \exp\left(   -\frac{l^2}{2v}  \right)\times 
 \frac{\lambda^2}{2} \exp\left(  - \frac{\lambda^2 v}{2} \right) \d v
=   \frac{\lambda}{2} \exp\left( -\lambda |l|  \right).
$$

\begin{proof}
By completing the square in the exponential term, we have
$$
f(l) = \int_0^\infty 
 \frac{\lambda^2 e^{-\lambda |l|}  }{2 \sqrt{2\pi}}  \frac{1}{\sqrt{v}} 
  \exp\left\{    - \frac{1}{2}   \left(   \frac{|l|}{\sqrt{v}}   - \lambda  \sqrt{v}  \right)^2      \right\} \d v.
$$
Change of variable $u=\sqrt{v}$ with $\d v = 2 u \d u$ gives
$$
f(l) = 
 \frac{\lambda^2 e^{-\lambda |l|}  }{  \sqrt{2\pi}}  
 \int_0^\infty 
  \exp\left\{    - \frac{1}{2}   \left(   \frac{|l|}{u}   - \lambda  u \right)^2      \right\} \d u.
$$
Change of variable 
$$
\eta = |l|/u - \lambda u, \quad 
\frac{\d u}{\d \eta} = \left(  - 1+ \frac{\eta}{  \sqrt{  \eta^2 + 4\lambda |l|      } }  \right) \Big / (2\lambda)
$$
gives
\begin{eqnarray*}
f(l) &=&
 \frac{\lambda^2 e^{-\lambda |l|}  }{  2\lambda  \sqrt{2\pi}}  
  \int_{-\infty }^\infty 
  \exp\left(     - \frac{\eta^2}{2}     \right)    \left(   1- \frac{\eta}{  \sqrt{  \eta^2 + 4\lambda |l|      } }  \right)     \d \eta \\
  &=& 
   \frac{\lambda e^{-\lambda |l|}  }{  2}    \int_{-\infty }^\infty    \frac{1}{\sqrt{2\pi}}   \exp\left(     - \frac{\eta^2}{2}     \right)  \d \eta
   -     \frac{\lambda e^{-\lambda |l|}  }{  2 \sqrt{2\pi}}    \int_{-\infty }^\infty  \exp\left(     - \frac{\eta^2}{2}     \right)    \frac{\eta}{  \sqrt{  \eta^2 + 4\lambda |l|      } } \d \eta
   =  \frac{\lambda e^{-\lambda |l|}  }{  2} ,
\end{eqnarray*}
where the last line follows from the standard Normal density and the integral of an odd function. 
\end{proof}

\end{document}